Correction of '*J. Laderman, V. Pan, X.-H. Sha, On practical Algorithms for Accelerated Matrix Multiplication, Linear Algebra and its Applications. Vol. 162-164 (1992) pp. 557-588*' [1]

Jerzy S. Respondek

*Silesian University of Technology, Faculty of Automatic Control, Electronics and Computer Science Applied Informatics Department, ul. Akademicka 16, 44-100 Gliwice, Poland, jrespondek@polsl.pl*

**Abstract**

In this article we corrected the trilinear formula for triple disjoint matrix multiplication given in the article '*J. Laderman, V. Pan, X. H. Sha, On practical Algorithms for Accelerated Matrix Multiplication, Linear Algebra and its Applications. Vol. 162-164 (1992) pp. 557-588*', which is incorrect for matrix dimensions equal to two or greater. That formula is a base of two algorithms, for disjoint and single matrix multiplication. The necessary correction made the amount of non-scalar products raise, slightly increasing the algorithm time complexity. We also corrected explicit formulas, in the bilinear form, for triple disjoint matrix multiplication. They are explicit, thus convenient for practical use of fast matrix multiplication algorithms in question.

## 1. Introduction

The article [4], as a base for practical algorithms for triple disjoint and single matrix multiplication, gave trilinear identity (1) with no proof. We show that it needs to be corrected to find correct values of certain constant parameters.

To find correct values of those parameters, and simultaneously correct algorithms in the bilinear form, we proposed proof of the identity (1). The article [4] proposed for the $f, g$ parameters to use particular values, which caused that a few terms vanished, this way improving the efficiency. We proved that in the correct form of trilinear identity (1) such optimization is not possible.

The paper is organized as follows: in section 2 we performed the proof of trilinear identity (1) finding correct values of a few parameters, in section 3 we presented the algorithm for triple disjoint matrix multiplication in the bilinear form, in section 4 we determined the complexities for triple disjoint and single matrix multiplication yielding from identity (1) and in section 5 we discussed the results.

---
[1] This work is supported by the National Research Fund No 02/100/BK_23/0027

## 2. Trilinear representation of three disjoint matrix products

Our objective is to calculate three disjoint matrix products $C = AB$, $W = UV$ and $Z = XY$, where $A = [a_{ij}]$, $B = [b_{jk}]$, $U = [u_{jk}]$, $V = [v_{ki}]$, $X = [x_{ki}]$, $Y = [y_{ij}]$, and $i, j, k$ range from 1 to $n$. The article [4], on pages 565-566 of section 3, proposes the following trilinear identity, with the $q$ parameter equal to 1, representing three disjoint matrix products[2]:

$$\mathrm{Trace}(ABC + UVW + XYZ) = \sum_{i,j,k} (a_{ij} b_{jk} c_{ki} + u_{jk} v_{ki} w_{ij} + x_{ki} y_{ij} z_{jk}) =$$

$$= \sum_{i,j,k} \Big[ (a_{ij} + u_{jk} + x_{ki})(b_{jk} + v_{ki} + y_{ij})(c_{ki} + w_{ij} + z_{jk}) +$$
$$\quad - T_0(i,j,k) - T_1(j,k,i) - T_2(k,i,j) \Big]$$
$$+ T_0' + T_1' + T_2' +$$

$$\left.\begin{aligned}
&-\sum_{i,j}\left(a_{ij} + \frac{1}{g}\sum_k x_{ki}\right)\left(y_{ij} + \frac{1}{g}\sum_k v_{ki}\right)\left(gw_{ij} + \sum_k c_{ki}\right) + \\
&-\sum_{i,j}\left(a_{ij} + \frac{1}{h}\sum_k u_{jk}\right)\left(y_{ij} + \frac{1}{h}\sum_k b_{jk}\right)\left(hw_{ij} + \sum_k z_{jk}\right) + \\
&-\sum_{j,k}\left(u_{jk} + \frac{1}{g}\sum_i a_{ij}\right)\left(b_{jk} + \frac{1}{g}\sum_i y_{ij}\right)\left(gz_{jk} + \sum_i w_{ij}\right) + \\
&-\sum_{j,k}\left(u_{jk} + \frac{1}{h}\sum_i x_{ki}\right)\left(b_{jk} + \frac{1}{h}\sum_i v_{ki}\right)\left(hz_{jk} + \sum_i c_{ki}\right) + \\
&-\sum_{k,i}\left(x_{ki} + \frac{1}{g}\sum_j u_{jk}\right)\left(v_{ki} + \frac{1}{g}\sum_j b_{jk}\right)\left(gc_{ki} + \sum_j z_{jk}\right) + \\
&-\sum_{k,i}\left(x_{ki} + \frac{1}{h}\sum_j a_{ij}\right)\left(v_{ki} + \frac{1}{h}\sum_j y_{ij}\right)\left(hc_{ki} + \sum_j w_{ij}\right)
\end{aligned}\right\} \quad (*) \quad\quad (1)$$

for $g + h = n$. The sub-expressions are the following:

$$T_0(i,j,k) = a_{ij} v_{ki} z_{jk}, \quad\quad T_1(j,k,i) = u_{jk} y_{ij} c_{ki}, \quad\quad T_2(k,i,j) = x_{ki} b_{jk} w_{ij}$$

$$T_0' = \sum_i \Bigg[ \left(\frac{1}{h}\sum_j a_{ij} + \frac{1}{g}\sum_k x_{ki}\right)\left(\frac{1}{h}\sum_j y_{ij} + \frac{1}{g}\sum_k v_{ki}\right)\left(g\sum_j w_{ij} + h\sum_k c_{ki}\right) +$$
$$+ \frac{q-g}{h^2}\sum_j a_{ij} \sum_j y_{ij} \sum_j w_{ij} + \frac{q-h}{g^2}\sum_k x_{ki} \sum_k v_{ki} \sum_k c_{ki} \Bigg]$$

---

[2] We changed the notation from $T_0(i)$, $T_1(j)$ and $T_2(k)$ to $T_0'$, $T_1'$ and $T_2'$, because these expressions do not depend on $i, j$ and $k$, respectively.

$$T_1' = \sum_j \left[ \left( \frac{1}{g} \sum_i a_{ij} + \frac{1}{h} \sum_k u_{jk} \right) \left( \frac{1}{g} \sum_i y_{ij} + \frac{1}{h} \sum_k b_{jk} \right) \left( h \sum_i w_{ij} + g \sum_k z_{jk} \right) + \right.$$

$$\left. + \frac{q-h}{g^2} \sum_i a_{ij} \sum_i y_{ij} \sum_i w_{ij} + \frac{q-g}{h^2} \sum_k u_{jk} \sum_k b_{jk} \sum_k z_{jk} \right]$$

$$T_2' = \sum_k \left[ \left( \frac{1}{h} \sum_i x_{ki} + \frac{1}{g} \sum_j u_{jk} \right) \left( \frac{1}{h} \sum_i v_{ki} + \frac{1}{g} \sum_j b_{jk} \right) \left( g \sum_i c_{ki} + h \sum_j z_{jk} \right) + \right.$$

$$\left. + \frac{q-g}{h^2} \sum_i x_{ki} \sum_i v_{ki} \sum_i c_{ki} + \frac{q-h}{g^2} \sum_j u_{jk} \sum_j b_{jk} \sum_j z_{jk} \right]$$

To find the value of $q$, we will carry out the proof of (1). Let us start from the algebraic identity (5) from section 5 of the article [6], page 171:

$$a_{ij} b_{jk} c_{ki} + u_{jk} v_{ki} w_{ij} + x_{ki} y_{ij} z_{jk} = (a_{ij} + u_{jk} + x_{ki})(b_{jk} + v_{ki} + y_{ij})(c_{ki} + w_{ij} + z_{jk}) +$$

$$- \left[ a_{ij} y_{ij} (c_{ki} + w_{ij} + z_{jk}) + u_{jk} b_{jk} (w_{ij} + z_{jk} + c_{ki}) + x_{ki} v_{ki} (z_{jk} + c_{ki} + w_{ij}) \right] +$$

$$- \left[ a_{ij} \left( b_{jk} + v_{ki} \right) w_{ij} + u_{jk} \left( v_{ki} + y_{ij} \right) z_{jk} + x_{ki} \left( y_{ij} + b_{jk} \right) c_{ki} + \right.$$

$$\left. + \left( a_{ij} + u_{jk} \right) v_{ki} c_{ki} + \left( u_{jk} + x_{ki} \right) y_{ij} w_{ij} + \left( x_{ki} + a_{ij} \right) b_{jk} z_{jk} \right] + \quad (2)$$

$$- \left[ a_{ij} v_{ki} z_{jk} + u_{jk} y_{ij} c_{ki} + x_{ki} b_{jk} w_{ij} \right]$$

Let us sum both sides of (2) with respect to $i, j$ and $k$. We arrive at formulas (6, 7) of [6]:

$$\text{Trace}(ABC + UVW + XYZ) = \sum_{i,j,k} (a_{ij} b_{jk} c_{ki} + u_{jk} v_{ki} w_{ij} + x_{ki} y_{ij} z_{jk}) =$$

$$= \sum_{i,j,k} (a_{ij} + u_{jk} + x_{ki})(b_{jk} + v_{ki} + y_{ij})(c_{ki} + w_{ij} + z_{jk}) + \quad (3)$$

$$\left. \begin{array}{l} - \sum_{i,j} a_{ij} y_{ij} \sum_k (c_{ki} + w_{ij} + z_{jk}) - \sum_{j,k} u_{jk} b_{jk} \sum_i (w_{ij} + z_{jk} + c_{ki}) + \\ - \sum_{k,i} x_{ki} v_{ki} \sum_j (z_{jk} + c_{ki} + w_{ij}) + \end{array} \right\} \text{(A)}$$

$$\left. \begin{array}{l} - \sum_{i,j} a_{ij} \left[ \sum_k \left( b_{jk} + v_{ki} \right) \right] w_{ij} - \sum_{j,k} u_{jk} \left[ \sum_i \left( v_{ki} + y_{ij} \right) \right] z_{jk} + \\ - \sum_{k,i} x_{ki} \left[ \sum_j \left( y_{ij} + b_{jk} \right) \right] c_{ki} - \sum_{k,i} \left[ \sum_j \left( a_{ij} + u_{jk} \right) \right] v_{ki} c_{ki} + \\ - \sum_{i,j} \left[ \sum_k \left( u_{jk} + x_{ki} \right) \right] y_{ij} w_{ij} - \sum_{j,k} \left[ \sum_i \left( x_{ki} + a_{ij} \right) \right] b_{jk} z_{jk} + \end{array} \right\} \text{(B)}$$

$$\left. - \sum_{i,j,k} a_{ij} v_{ki} z_{jk} - \sum_{i,j,k} u_{jk} y_{ij} c_{ki} - \sum_{i,j,k} x_{ki} b_{jk} w_{ij} \right\} \text{(C)}$$

We will show that the trilinear formula (1) reduces to identity (3). For this purpose we divide factors in (1) into separate, complementary groups.

- As the first group of factors let us have a look at a series of sums indicated by (*) of (1), the following three sums of triple products:

$$\sum_{i,j} a_{ij} y_{ij} \left( g w_{ij} + \sum_{k} c_{ki} \right) + \sum_{i,j} a_{ij} y_{ij} \left( h w_{ij} + \sum_{k} z_{jk} \right)$$

$$\sum_{j,k} u_{jk} b_{jk} \left( g z_{jk} + \sum_{i} w_{ij} \right) + \sum_{j,k} u_{jk} b_{jk} \left( h z_{jk} + \sum_{i} c_{ki} \right) \quad (4)$$

$$\sum_{k,i} x_{ki} v_{ki} \left( g c_{ki} + \sum_{j} z_{jk} \right) + \sum_{k,i} x_{ki} v_{ki} \left( h c_{ki} + \sum_{j} w_{ij} \right)$$

Considering $g + h = n$ we obtain, by the example of the first sum:

$$\sum_{i,j} a_{ij} y_{ij} \left( g w_{ij} + \sum_{k} c_{ki} \right) + \sum_{i,j} a_{ij} y_{ij} \left( h w_{ij} + \sum_{k} z_{jk} \right) =$$

$$= \sum_{i,j} a_{ij} y_{ij} \left( g w_{ij} + h w_{ij} + \sum_{k} c_{ki} + \sum_{k} z_{jk} \right) = \sum_{i,j} a_{ij} y_{ij} \left( n w_{ij} + \sum_{k} c_{ki} + \sum_{k} z_{jk} \right)$$

We can rewrite $n w_{ij}$ in the form $\sum_{k} w_{ij}$. So:

$$\sum_{i,j} a_{ij} y_{ij} \left( n w_{ij} + \sum_{k} c_{ki} + \sum_{k} z_{jk} \right) = \sum_{i,j} a_{ij} y_{ij} \sum_{k} (c_{ki} + w_{ij} + z_{jk})$$

Thus this part of proven equality (1) is reduced to the first sum of part (A) in (3). The remaining two expressions in (4) analogically represent the second and third sum of part (A) in (3).

- For the second group of terms we again look at from the series (*) in (1) six triple product sums of the following form:

$$\sum_{i,j} \left( \frac{1}{g} \sum_{k} x_{ki} \right) \left( \frac{1}{g} \sum_{k} v_{ki} \right) g w_{ij}, \qquad \sum_{i,j} \left( \frac{1}{h} \sum_{k} u_{jk} \right) \left( \frac{1}{h} \sum_{k} b_{jk} \right) h w_{ij}$$

$$\sum_{j,k} \left( \frac{1}{g} \sum_{i} a_{ij} \right) \left( \frac{1}{g} \sum_{i} y_{ij} \right) g z_{jk}, \qquad \sum_{j,k} \left( \frac{1}{h} \sum_{i} x_{ki} \right) \left( \frac{1}{h} \sum_{i} v_{ki} \right) h z_{jk} \quad (5)$$

$$\sum_{k,i} \left( \frac{1}{g} \sum_{j} u_{jk} \right) \left( \frac{1}{g} \sum_{j} b_{jk} \right) g c_{ki}, \qquad \sum_{k,i} \left( \frac{1}{h} \sum_{j} a_{ij} \right) \left( \frac{1}{h} \sum_{j} y_{ij} \right) h c_{ki}$$

To get to know which of the above terms (5) corresponds to the terms in (3) and/or is reduced by the correction terms $T_0', T_1', T_2'$ and $T_0(i,j,k)$, $T_1(j,k,i)$, $T_2(k,i,j)$, we look at the example of a selected expression from (5):

$$\sum_{k,i}\left(\frac{1}{h}\sum_j a_{ij}\right)\left(\frac{1}{h}\sum_j y_{ij}\right)hc_{ki} = \sum_i\left(\frac{1}{h}\sum_j a_{ij}\right)\left(\frac{1}{h}\sum_j y_{ij}\right)h\sum_k c_{ki}$$

It is now clear that each of the expressions (5) is reduced by its matching part in the terms $T_0', T_1'$ or $T_2'$ just in the trilinear equality (1).

- Another group of terms includes expressions from the series (*) in (1) of the form:

$$\sum_{i,j} a_{ij}\left(\frac{1}{h}\sum_k b_{jk}\right)\left(hw_{ij}+\sum_k z_{jk}\right) + \sum_{j,k}\left(\frac{1}{g}\sum_i a_{ij}\right)b_{jk}\left(gz_{jk}+\sum_i w_{ij}\right)$$

$$\sum_{j,k} u_{jk}\left(\frac{1}{h}\sum_i v_{ki}\right)\left(hz_{jk}+\sum_i c_{ki}\right) + \sum_{k,i}\left(\frac{1}{g}\sum_j u_{jk}\right)v_{ki}\left(gc_{ki}+\sum_j z_{jk}\right)$$

$$\sum_{k,i} x_{ki}\left(\frac{1}{g}\sum_j b_{jk}\right)\left(gc_{ki}+\sum_j z_{jk}\right) + \sum_{j,k}\left(\frac{1}{h}\sum_i x_{ki}\right)b_{jk}\left(hz_{jk}+\sum_i c_{ki}\right) \quad (6)$$

$$\sum_{i,j} a_{ij}\left(\frac{1}{g}\sum_k v_{ki}\right)\left(gw_{ij}+\sum_k c_{ki}\right) + \sum_{k,i}\left(\frac{1}{h}\sum_j a_{ij}\right)v_{ki}\left(hc_{ki}+\sum_j w_{ij}\right)$$

$$\sum_{j,k} u_{jk}\left(\frac{1}{g}\sum_i y_{ij}\right)\left(gz_{jk}+\sum_i w_{ij}\right) + \sum_{i,j}\left(\frac{1}{h}\sum_k u_{jk}\right)y_{ij}\left(hw_{ij}+\sum_k z_{jk}\right)$$

$$\sum_{k,i} x_{ki}\left(\frac{1}{h}\sum_j y_{ij}\right)\left(hc_{ki}+\sum_j w_{ij}\right) + \sum_{i,j}\left(\frac{1}{g}\sum_k x_{ki}\right)y_{ij}\left(gw_{ij}+\sum_k c_{ki}\right)$$

By the example of the first expression in series (6) we give a general rule:

$$\sum_{i,j} a_{ij}\left(\frac{1}{h}\sum_k b_{jk}\right)\left(hw_{ij}+\sum_k z_{jk}\right) + \sum_{j,k}\left(\frac{1}{g}\sum_i a_{ij}\right)b_{jk}\left(gz_{jk}+\sum_i w_{ij}\right) =$$

$$= \sum_{i,j,k}\left[a_{ij}b_{jk}\left(w_{ij}+\frac{1}{h}\sum_{k_1} z_{jk_1}\right) + a_{ij}b_{jk}\left(z_{jk}+\frac{1}{g}\sum_{i_1} w_{i_1 j}\right)\right] =$$

$$= \sum_{i,j,k}\left[\left(a_{ij}b_{jk}w_{ij}+a_{ij}b_{jk}z_{jk}\right) + a_{ij}b_{jk}\left(\frac{1}{h}\sum_{k_1} z_{jk_1}+\frac{1}{g}\sum_{i_1} w_{i_1 j}\right)\right]$$

We can see that each of the six sums (6) produces two groups of expressions:

- Sum of the expressions of the form $a_{ij}b_{jk}w_{ij}+a_{ij}b_{jk}z_{jk}$. It is easy to see that they are desired terms, which summed over all six expressions (6) forms part (B) of equality (3).
- Sum of the expressions of the form $a_{ij}b_{jk}\left(1/h\sum_{k_1} z_{jk_1}+1/g\sum_{i_1} w_{i_1 j}\right)$ which are undesired. They need to be reduced by correction terms of trilinear equality (1). Let us find a corresponding expression in (1), giving terms of a general form $\sum a_{ij}\sum b_{jk}\sum w_{ij}$ and $\sum a_{ij}\sum b_{jk}\sum z_{jk}$.

We can find a proper expression in the correction term $T_1'$:

$$\sum_j \left(\frac{1}{g}\sum_i a_{ij}\right)\left(\frac{1}{h}\sum_k b_{jk}\right)\left(h\sum_i w_{ij} + g\sum_k z_{jk}\right) = \sum_j \sum_i \sum_k \frac{1}{gh} a_{ij} b_{jk} \left(h\sum_{i_1} w_{i_1 j} + g\sum_{k_1} z_{jk_1}\right) =$$

$$= \sum_{i,j,k} a_{ij} b_{jk} \left(\frac{1}{h}\sum_{k_1} z_{jk_1} + \frac{1}{g}\sum_{i_1} w_{i_1 j}\right)$$

To sum up, the group of six terms of the form (6) create group (B) of equality (3) and the undesired expressions they produce are corrected by the corresponding correction terms of $T_0', T_1', T_2'$.

It is also worth to notice that all the terms in (3) are covered by the so-far presented groups of expressions from trilinear identity (1). It means that the remaining terms in (1) reduce themselves within that identity.

- The last group of terms in the series (*) in (1) comprises six sums of triple products, of a common form:

$$\sum_{i,j}\left(\frac{1}{g}\sum_k x_{ki}\right)\left(\frac{1}{g}\sum_k v_{ki}\right)\left(\sum_k c_{ki}\right), \qquad \sum_{i,j}\left(\frac{1}{h}\sum_k u_{jk}\right)\left(\frac{1}{h}\sum_k b_{jk}\right)\left(\sum_k z_{jk}\right)$$

$$\sum_{j,k}\left(\frac{1}{g}\sum_i a_{ij}\right)\left(\frac{1}{g}\sum_i y_{ij}\right)\left(\sum_i w_{ij}\right), \qquad \sum_{j,k}\left(\frac{1}{h}\sum_i x_{ki}\right)\left(\frac{1}{h}\sum_i v_{ki}\right)\left(\sum_i c_{ki}\right) \qquad (7)$$

$$\sum_{k,i}\left(\frac{1}{g}\sum_j u_{jk}\right)\left(\frac{1}{g}\sum_j b_{jk}\right)\left(\sum_j z_{jk}\right), \qquad \sum_{k,i}\left(\frac{1}{h}\sum_j a_{ij}\right)\left(\frac{1}{h}\sum_j y_{ij}\right)\left(\sum_j w_{ij}\right)$$

Let us analyze the first trilinear product in (7):

$$\sum_{i,j}\left(\frac{1}{g}\sum_k x_{ki}\right)\left(\frac{1}{g}\sum_k v_{ki}\right)\left(\sum_k c_{ki}\right)$$

The inside, summed expression $\left(1/g\sum_k x_{ki}\right)\left(1/g\sum_k v_{ki}\right)\left(\sum_k c_{ki}\right)$ does not depend on the $j$ variable. So we can claim:

$$\sum_{i,j}\left(\frac{1}{g}\sum_k x_{ki}\right)\left(\frac{1}{g}\sum_k v_{ki}\right)\left(\sum_k c_{ki}\right) = n \cdot \sum_i \left(\frac{1}{g}\sum_k x_{ki}\right)\left(\frac{1}{g}\sum_k v_{ki}\right)\left(\sum_k c_{ki}\right) \qquad (8)$$

Thus the corresponding term in $T_0'$ reduces (8) iff $q = n$. For the remaining sums of triple products (7) there exist corresponding five reducing expressions in $T_0', T_1'$ and $T_2'$.

- The last group of terms are in $T_0', T_1'$ and $T_2'$. These terms balance themselves in the given term, in the outside sums. The terms in question are the following:

$$T_0': \frac{g}{h^2}\sum_j a_{ij} \sum_j y_{ij} \sum_j w_{ij}, \quad \frac{h}{g^2}\sum_k x_{ki} \sum_k v_{ki} \sum_k c_{ki}$$

$$T_1': \frac{h}{g^2}\sum_i a_{ij} \sum_i y_{ij} \sum_i w_{ij}, \quad \frac{g}{h^2}\sum_k u_{jk} \sum_k b_{jk} \sum_k z_{jk}$$

$$T_2': \frac{g}{h^2}\sum_i x_{ki} \sum_i v_{ki} \sum_i c_{ki}, \quad \frac{h}{g^2}\sum_j u_{jk} \sum_j b_{jk} \sum_j z_{jk}$$

To sum up, in the trilinear identity in section 3 of the article [4], pages 565-566, its part included in the correction terms $T_0', T_1'$ and $T_2'$ of a general form:

$$\frac{1-g}{h^2}\sum f^{(1)}_{i_1 i_2} \sum f^{(2)}_{i_1 i_2} \sum f^{(3)}_{i_1 i_2}, \quad \text{and} \quad \frac{1-h}{g^2}\sum f^{(4)}_{i_1 i_2} \sum f^{(5)}_{i_1 i_2} \sum f^{(6)}_{i_1 i_2}$$

requires to be corrected to, respectively:

$$\frac{n-g}{h^2}\sum f^{(1)}_{i_1 i_2} \sum f^{(2)}_{i_1 i_2} \sum f^{(3)}_{i_1 i_2}, \quad \text{and} \quad \frac{n-h}{g^2}\sum f^{(4)}_{i_1 i_2} \sum f^{(5)}_{i_1 i_2} \sum f^{(6)}_{i_1 i_2}$$

### 3. Bilinear form of the algorithm

In this section we show the bilinear form of the algorithm which is ready for practical purposes. It is obtained by equating the coefficients of $c_{ki}, w_{ij}$ and $z_{jk}$ on both sides. Like in [4] we assume that $g = 1$ and $h = n - 1$ as well as we omit the $T_0(i,j,k)$, $T_1(j,k,i)$ and $T_2(k,i,j)$ terms.

At first let us define the non-scalar products. We have to perform $n^3$ of $P_{ijk}$ products:

$$P_{ijk} = (a_{ij} + u_{jk} + x_{ki})(b_{jk} + v_{ki} + y_{ij})$$

Next come six series with $n^2$ products in each of them, of the form:

$$P^{(0)}_{ij} = \sum_{i,j}\left(a_{ij} + \sum_k x_{ki}\right)\left(y_{ij} + \sum_k v_{ki}\right), \quad P^{(1)}_{ij} = \sum_{i,j}\left(a_{ij} + \frac{1}{n-1}\sum_k u_{jk}\right)\left(y_{ij} + \frac{1}{n-1}\sum_k b_{jk}\right)$$

$$P^{(2)}_{jk} = \sum_{j,k}\left(u_{jk} + \sum_i a_{ij}\right)\left(b_{jk} + \sum_i y_{ij}\right), \quad P^{(3)}_{jk} = \sum_{j,k}\left(u_{jk} + \frac{1}{n-1}\sum_i x_{ki}\right)\left(b_{jk} + \frac{1}{n-1}\sum_i v_{ki}\right)$$

$$P^{(4)}_{ki} = \sum_{k,i}\left(x_{ki} + \sum_j u_{jk}\right)\left(v_{ki} + \sum_j b_{jk}\right), \quad P^{(5)}_{ki} = \sum_{k,i}\left(x_{ki} + \frac{1}{n-1}\sum_j a_{ij}\right)\left(v_{ki} + \frac{1}{n-1}\sum_j y_{ij}\right)$$

Finally from $T_0'$, $T_1'$ and $T_2'$ we have nine groups with $n$ products in each, of the form:

$$P_i^{(0)} = \left(\frac{1}{n-1}\sum_j a_{ij} + \sum_k x_{ki}\right)\left(\frac{1}{n-1}\sum_j y_{ij} + \sum_k v_{ki}\right)$$

$$P_i^{(1)} = \frac{1}{n-1}\sum_j a_{ij}\sum_j y_{ij}\sum_j w_{ij}, \quad P_i^{(2)} = \sum_k x_{ki}\sum_k v_{ki}\sum_k c_{ki}$$

$$P_j^{(3)} = \left(\sum_i a_{ij} + \frac{1}{n-1}\sum_k u_{jk}\right)\left(\sum_i y_{ij} + \frac{1}{n-1}\sum_k b_{jk}\right)$$

$$P_j^{(4)} = \sum_i a_{ij}\sum_i y_{ij}\sum_i w_{ij}, \quad P_j^{(5)} = \frac{1}{n-1}\sum_k u_{jk}\sum_k b_{jk}\sum_k z_{jk}$$

$$P_k^{(6)} = \left(\frac{1}{n-1}\sum_i x_{ki} + \sum_j u_{jk}\right)\left(\frac{1}{n-1}\sum_i v_{ki} + \sum_j b_{jk}\right)$$

$$P_k^{(7)} = \frac{1}{n-1}\sum_i x_{ki}\sum_i v_{ki}\sum_i c_{ki}, \quad P_k^{(8)} = \sum_j u_{jk}\sum_j b_{jk}\sum_j z_{jk}$$

The desired three matrix products are expressed explicitly by linear combinations of the above products:

$$(C = AB)_{ik} = \sum_j \left(P_{ijk} - P_{ij}^{(0)} - P_{jk}^{(3)}\right) - \left[P_{ki}^{(4)} + (n-1)P_{ki}^{(5)}\right] + (n-1)P_i^{(0)} + P_i^{(2)} + P_k^{(6)} + P_k^{(7)}$$

$$(W = UV)_{ji} = \sum_k \left(P_{ijk} - P_{jk}^{(2)} - P_{ki}^{(5)}\right) - \left[P_{ij}^{(0)} + (n-1)P_{ij}^{(1)}\right] + P_i^{(0)} + P_i^{(1)} + (n-1)P_j^{(3)} + P_j^{(4)}$$

$$(Z = XY)_{kj} = \sum_i \left(P_{ijk} - P_{ij}^{(1)} - P_{ki}^{(4)}\right) - \left[P_{jk}^{(2)} + (n-1)P_{jk}^{(3)}\right] + P_j^{(3)} + P_j^{(5)} + (n-1)P_k^{(6)} + P_k^{(8)}$$

## 4. Time complexity issues

Without the $T_0(i,j,k)$, $T_1(j,k,i)$ and $T_2(k,i,j)$ terms, trilinear identity (1) contains $n^3 + 6n^2 + 9n$ products. Taking like in the article [4] $g = 1$ and $h = n-1$, or any other non-zero values satisfying $g+h = n$, does not cause $3n$ terms to vanish. Applying the $2 \times 2$ block scheme for $n \times n$ matrices from section 4 of the article [4] we get rid of the $T_0(i,j,k)$, $T_1(j,k,i)$ and $T_2(k,i,j)$ terms and obtain:

$$M_{\text{disjoint}}(2n) = 8(n^3 + 6n^2 + 9n) = (2n)^3 + 12(2n)^2 + 36(2n)$$

products necessary to calculate. The obtained exponent expression is $\omega \leq \log_{2n}\left[M_{\text{disjoint}}(2n)/3\right]$. For $2n = 48$ we obtain an algorithm with exponent $\omega < 2.77706$

Applying trilinear identity (1) for single matrix multiplication, in a way from section 6 of the article [4], we have to calculate $(n^3 - n)/3 + 2n^2 + 3n$ products. The $2 \times 2$ block scheme for $n \times n$ matrices leads here to

$$M_{\text{single}}(2n) = \frac{8}{3}(n^3 - n) + 16n^2 + 24n = \frac{(2n)^3}{3} + 4(2n)^2 + \frac{32}{3}(2n)$$

terms. The exponent inequality here is $\omega \leq \log_{2n}\left[M_{\text{single}}(2n)\right]$ with minimum also for $2n = 48$ equal to $\omega < 2.776706$.

## 5. Conclusions and perspectives

Since the publication of Strassen's work [7], which gave the first non-commutative algorithm for matrix multiplication with time complexity exponent smaller than definitional 3, the exponent was improved a series of times currently reaching $\omega = 2.3728596$ in the work [1]. The exponents yielding from the algorithms proposed in this article and the article [4] are slightly above $2.77$, but they are easier to use.

As an example of the technique used in algorithms giving better exponents we can mention the so-called approximate algorithms, devised by Bini, Capovani et al. [2]. The approximate algorithms can be transformed into exact ones by the technique presented by Bini [3], but it induces an additional logarithmic term in the obtained complexity. Another technique is to apply disjoint matrix multiplication for a single matrix multiplication, but to be effective it requires a series of recursion steps. The general theorem to determine the complexity exponent by this technique, as well as for the combination with approximate methods, is given in the article [5].

On the opposite, algorithms presented in [4] corrected by this article require neither recursion nor approximate techniques to reduce the number of non-scalar multiplication to be performed. Thus they have potential to be convenient in practice.

**Acknowledgement**

I would like to thank my university colleagues for stimulating discussions.